\documentclass{amsart}
\theoremstyle{plain}
\newtheorem{thm}{Theorem}[section]
\newtheorem{lem}[thm]{Lemma}
\theoremstyle{definition}
\newtheorem{dfn}[thm]{Definition}
\newtheorem{exa}[thm]{Example}
\newtheorem{rem}[thm]{Remark}

\newcommand{\N}{\mathbb{N}}

\DeclareMathOperator{\T}{T}
\DeclareMathOperator{\GL}{GL}
\DeclareMathOperator{\tr}{tr}

\DeclareMathOperator{\Rep}{Rep}
\DeclareMathOperator{\Der}{Der}
\DeclareMathOperator{\CDer}{CDer}
\DeclareMathOperator{\Inn}{Inn}
\DeclareMathOperator{\Mod}{Mod}

\DeclareMathOperator{\dbslash}{/\!\!/}

\newcommand{\kdelta}{\mbox{\boldmath $\delta$}}

\newcommand{\lp}{\langle}
\newcommand{\rp}{\rangle}

\begin{document}
\title{A note on noncommutative Poisson structures}

\author{William Crawley-Boevey}
\address{Department of Pure Mathematics, University of Leeds, Leeds LS2 9JT, UK}
\email{w.crawley-boevey@leeds.ac.uk}

\thanks{Mathematics Subject Classification (2000): Primary 16W25, 17B63, 53D30}

\begin{abstract}
We introduce a new type of noncommutative Poisson structure on associative
algebras. It induces Poisson structures on the moduli spaces classifying
semisimple modules. Path algebras of doubled quivers and preprojective algebras
have noncommutative Poisson structures given by the necklace Lie algebra.
\end{abstract}
\maketitle

Recall that a Poisson bracket on a commutative ring $A$ is a Lie bracket
\[
\{ -,-\} : A\times A\to A
\]
which satisfies the Leibnitz rule
\[
\{ a,bc \} = b\{a,c\} + \{a,b\}c.
\]
The same definition can be used for a Poisson bracket on
a noncommutative ring, but it seems that this is too restrictive,
as by a theorem of Farkas and Letzter \cite{FL} the only Poisson
brackets on a genuinely noncommutative prime ring are the commutator
bracket $[a,b]=ab-ba$ and multiples of it (in a suitable sense).

A notion of a noncommutative Poisson structure has
been suggested by Xu~\cite{Xu} and Block and Getzler~\cite{BG}.
It has the property that if $A$ has a noncommutative Poisson
structure, then the centre $Z(A)$ has a Poisson bracket, but
otherwise the relationship with Poisson brackets is unclear.

In this paper we introduce a new type of noncommutative Poisson structure.
It is the weakest structure we can find which (when $A$ is a finitely generated
$K$-algebra and $K$ is an algebraically closed field of characteristic zero)
induces Poisson brackets on the coordinate rings of the moduli spaces
$\Mod(A,n)\dbslash \GL_n(K)$ classifying $n$-dimensional semisimple $A$-modules.
With this notion, path algebras of doubled quivers,
preprojective algebras and multiplicative preprojective
algebras all have noncommutative Poisson structures.

For a much deeper approach, see the work of Van den Bergh \cite{V}.
I would like to thank M.~Van den Bergh, who raised this problem,
and G.~Van de Weyer and Pu~Zhang for some useful
discussions.

\section{Definition and Examples}
Throughout, we work over a commutative base ring $K$ and,
where appropriate, maps are assumed to be $K$-linear.
Let $A$ be an associative $K$-algebra (with~1).
Recall that the zeroth Hochschild homology of $A$ is $A/[A,A]$,
where $[A,A]$ is the subset of $A$ spanned by the commutators.
We write $\overline a$ for the image of $a\in A$ in $A/[A,A]$.
Observe that if $d:A\to A$ is a derivation, then since
\[
d([a,b])= [a,d(b)] + [d(a),b] \in [A,A],
\]
there is an induced linear map $\overline{d}:A/[A,A]\to A/[A,A]$,
$\overline a \mapsto \overline{d(a)}$.

\begin{dfn}
By a \emph{noncommutative Poisson structure} on $A$ we mean a
Lie bracket $\lp-,-\rp$ on $A/[A,A]$, such that for each $a\in A$
the map
\[
\lp \overline a,-\rp:A/[A,A]\to A/[A,A]
\]
is induced by a derivation $d_a:A\to A$.
\end{dfn}

\begin{rem}
Observe that $\overline{d} = 0$ if $d$ is an inner derivation,
so if the first Hochschild cohomology of $A$ vanishes, any
noncommutative Poisson structure on $A$ must be zero.
\end{rem}

\begin{rem}
Writing $\Der(A,A)$ for the set of derivations $A\to A$,
and $\CDer(A,A)$ for the subset consisting of the derivations with
image contained in $[A,A]$ (which of course contains the set
$\Inn(A,A)$ of inner derivations), a noncommutative Poisson
structure determines a linear map
\[
\pi : A/[A,A] \to \Der(A,A)/\CDer(A,A)
\]
with $\pi(a)$ represented by the derivation $d_a$.
Conversely, a noncommutative Poisson structure may be defined
as a linear map $\pi$ such that the bracket $\lp-,-\rp$
on $A/[A,A]$ defined by $\lp \overline a,\overline b\rp = \overline{d_a(b)}$,
where $d_a$ is a representative of $\pi(a)$, is a Lie bracket.
\end{rem}

\begin{exa}
If $A$ is commutative, then a noncommutative Poisson structure on
$A$ is exactly the same as a Poisson bracket.
\end{exa}

\begin{exa}
Any Poisson bracket $\{-,-\}$ on $A$ induces a noncommutative Poisson
structure via $\lp \overline a,\overline b\rp = \overline{\{a,b\}}$.
\end{exa}

\begin{exa}
The necklace Lie algebra of \cite{BL} and \cite{G} is a noncommutative
Poisson structure. Let $Q$ be a quiver with vertex set $\{1,2,\dots,m\}$,
and let $\overline{Q}$ be its double, obtained by adjoining a reverse
arrow $a^*:w\to v$ for each arrow $a:v\to w$ in $Q$.
We extend $*$ to an involution on the set of all arrows in $\overline Q$
by defining $a^{**}=a$. Define $\epsilon(a)$ for all arrows $a\in\overline Q$
by $\epsilon(a)=1$ if $a\in Q$, and $\epsilon(a)=-1$ if $a^*\in Q$.
If $p$ is a path in $\overline{Q}$, we write $\ell(p)$ for its length.
For $1\le i\le \ell(p)$, we write $p=p_{<i}\, p_i\, p_{>i}$, where $p_i$
is an arrow and $p_{<i}$ and $p_{>i}$ are paths of lengths $i-1$ and
$\ell(p)-i$ respectively.

If $p$ is a path, then the assignment
\[
d_p(q) = \sum_{i=1}^{\ell(q)} \sum_{j=1}^{\ell(p)} \epsilon(q_i)
\kdelta_{q_i^*,p_j} q_{<i}\, p_{>j}\, p_{<j}\, q_{>i},
\]
for $q$ a path, defines a derivation $d_p:K\overline{Q}\to K\overline{Q}$.
Here $\kdelta$ is the Kronecker delta function, so
\[
\kdelta_{q_i^*,p_j}=\begin{cases}
1 & \text{(if $q_i^*=p_j$)} \\
0 & \text{(otherwise).}
\end{cases}
\]
The necklace Lie algebra on $K\overline{Q}/[K\overline{Q},K\overline{Q}]$
is given by the bracket
\[
\lp \overline p,\overline q\rp = \overline{d_p(q)}.
\]
This construction shows that it is a noncommutative Poisson structure.
\end{exa}

\begin{exa}
The deformed preprojective algebra \cite{CBH} of weight $\lambda\in K^m$
is the algebra
$\Pi^\lambda = K\overline{Q}/J$, where $J$ is the ideal generated
by $w-\lambda$,
\[
w = \sum_{a\in Q} [a,a^*] = \sum_{a\in\overline{Q}} \epsilon(a) aa^*
\]
and $\lambda$ is identified with the corresponding linear combination
$\sum_{v=1}^m\lambda_v e_v$ of trivial paths.

It is shown in \cite{CBEG} that the necklace Lie algebra structure
on $K\overline{Q}/[K\overline{Q},K\overline{Q}]$ descends to
a Lie algebra structure on $\Pi^\lambda/[\Pi^\lambda,\Pi^\lambda]$.
To show that this is a noncommutative Poisson structure, it suffices
to show that the derivation $d_p$ descends to a derivation
of $\Pi^\lambda$. Now
\[
\begin{split}
d_p(w-\lambda) &= \sum_{a\in \overline{Q}} \epsilon(a) (d_p(a)a^* + a d_p(a^*))
\\
&= \sum_{a\in\overline Q} \epsilon(a) \left(
\epsilon(a) \sum_{j=1}^{\ell(p)} \kdelta_{a^*,p_j} \, p_{>j} \, p_{<j} \, a^*
+ \epsilon(a^*) \sum_{j=1}^{\ell(p)} a \, \kdelta_{a,p_j} \, p_{>j} \, p_{<j}
\right)
\\
&= \sum_{j=1}^{\ell(p)} p_{>j} \, p_{<j} \, p_j
- \sum_{j=1}^{\ell(p)} p_j \, p_{>j} \, p_{<j} = 0.
\end{split}
\]
It follows that $d_p(J) \subseteq J$, as required.

\end{exa}

For more about these last two examples, and for multiplicative
preprojective algebras, see \cite{V}.

\section{Representation Schemes}
Let $A$ be a $K$-algebra and
suppose that $e_1,\dots,e_m$ is a complete set of
orthogonal idempotents in $A$, that is,
\[
e_i^2 = e_i,
\quad
e_ie_j = 0\ (i\neq j),
\quad
e_1+\dots+e_m=1.
\]
If $\alpha\in\N^m$, then there is an affine scheme $\Rep(A,\alpha)$
whose set of $S$-valued points, where $S$ is a commutative $K$-algebra,
is the set of $A\otimes_K S$-module structures on
\[
S^{\alpha_1}\oplus \dots \oplus S^{\alpha_m}
\]
such that the idempotents $e_v$ act as projection onto the $v$th summand.
Equivalently, it is the set of $K$-algebra homomorphisms $A\to M_n(S)$,
where $n=\alpha_1+\dots+\alpha_m$, such that the image of $e_v$
is the matrix $\Delta^v$, where
\[
\Delta^v_{ij} = \begin{cases}
1 & \text{(if $i=j$ and $\sum_{w<v} \alpha_w < i \le \sum_{w\le v} \alpha_w$)}
\\
0 & \text{(otherwise).}
\end{cases}
\]
The coordinate ring $K[\Rep(A,\alpha)]$ is generated by
elements $a_{ij}$ with $a\in A$ and $1\le i,j\le n$, where
subject to the relations
\[
(ab)_{ij} = \sum_{k=1}^n a_{ik} b_{kj},
\quad
(\lambda a+\mu b)_{ij} = \lambda a_{ij} + \mu b_{ij},
\quad
(e_v)_{ij} = \Delta^v_{ij}
\]
for all $a,b\in A$, $\lambda,\mu\in K$, $1\le i,j\le n$ and $1\le v\le m$.

\begin{rem}
As a special case
one gets the scheme $\Mod(A,n)$ of $n$-dimensional $A$-module
structures, for any $K$-algebra $A$, by taking $m=1$ and $e_1=1$.
\end{rem}

\begin{dfn}
For $a\in A$, we define the element
\[
\tr_{\alpha} a = \sum_{i=1}^n a_{ii} \in K[\Rep(A,\alpha)].
\]
Observe that $\tr_{\alpha}(ab) = \tr_{\alpha}(ba)$, so $\tr_{\alpha}$ induces
a map $A/[A,A]\to K[\Rep(A,\alpha)]$ which we also denote $\tr_{\alpha}$.
We define $\T(A,\alpha)$ to be the subalgebra of $K[\Rep(A,\alpha)]$
generated by the elements $\tr_{\alpha} a$.
\end{dfn}

\begin{rem}
Suppose that $K$ is an algebraically closed field and $A$ is a finitely
generated $K$-algebra. The group
\[
\GL(\alpha) = \GL_{\alpha_1}(K)\times\dots\times\GL_{\alpha_m}(K),
\]
embedded as diagonal blocks in $\GL_n(K)$, acts naturally on $\Rep(A,\alpha)$,
and the closed points of the affine quotient scheme
$\Rep(A,\alpha)\dbslash\GL(\alpha)$ classify isomorphism classes of
semisimple $A$-modules $M$ of dimension vector $\alpha$
(that is, with $\dim e_v M = \alpha_v$ for all $v$).

Now if $K$ has characteristic zero, $\T(A,\alpha)$
is the coordinate ring of this quotient scheme.
Namely, the coordinate ring is the ring of invariants
$K[\Rep(A,\alpha)]^{\GL(\alpha)}$, where the action of
$\GL(\alpha)$ is given by
\[
g.a_{ij} = \sum_{k=1}^n \sum_{\ell=1}^n g_{ik} (g^{-1})_{\ell j} a_{k\ell}
\quad
(g\in\GL(\alpha), a\in A, 1\le i,j\le n).
\]
Clearly $\tr_{\alpha} a$ is an invariant, and in fact the
elements $\tr_{\alpha} a$ generate the ring of invariants.
For a path algebra $KQ$ this holds by a theorem of
Le~Bruyn and Procesi~\cite{LP}. In general, there is a
surjective homomorphism from a path algebra $\theta:KQ\to A$,
and hence a surjective homomorphism
\[
\theta_\alpha : K[\Rep(KQ,\alpha)]\to K[\Rep(A,\alpha)]
\]
with $\theta_\alpha(x_{ij})=\theta(x)_{ij}$.
Since $\GL(\alpha)$ is linearly reductive, it induces a surjection
\[
K[\Rep(KQ,\alpha)]^{\GL(\alpha)}\to K[\Rep(A,\alpha)]^{\GL(\alpha)}.
\]
Thus $K[\Rep(A,\alpha)]^{\GL(\alpha)}$ is generated by the elements
$\theta_\alpha(\tr_{\alpha} x) = \tr_{\alpha} \theta(x)$,
and so it is equal to $\T(A,\alpha)$.
\end{rem}

\begin{lem}
Given a derivation $d:A \to A$ with $d(e_v)=0$ for all $v$,
there is a unique derivation $d_\alpha:K[\Rep(A,\alpha)] \to K[\Rep(A,\alpha)]$
with $d_\alpha(a_{ij}) = d(a)_{ij}$ for all $i,j$ and $a\in A$.
\end{lem}

\begin{proof}
We define $d_\alpha$ on the elements $a_{ij}$ by the indicated formula,
and it extends uniquely to a derivation on the polynomial ring $K[a_{ij}]$.
This descends to a derivation on $K[\Rep(A,\alpha)]$ since
\[
\begin{split}
d_\alpha((ab)_{ij})
&= d(ab)_{ij}
= (ad(b)+d(a)b)_{ij}
= \sum_{k=1}^n a_{ik}d(b)_{kj} + d(a)_{ik}b_{kj}
\\
&= \sum_{k=1}^n a_{ik}d_\alpha(b_{kj}) + d_\alpha(a_{ik})b_{kj}
= d_\alpha \biggl( \sum_{k=1}^n  a_{ik} b_{kj} \biggr),
\end{split}
\]
and clearly $d_\alpha((\lambda a+\mu b)_{ij}) = d_\alpha(\lambda a_{ij} + \mu b_{ij})$,
and $d_\alpha((e_v)_{ij}) = 0 = d_\alpha(\Delta^v_{ij})$.
\end{proof}

\begin{thm}
If $\lp-,-\rp$ is a noncommutative Poisson structure on $A$, then for
any $\alpha$ there is a unique Poisson bracket $\{-,-\}$ on $\T(A,\alpha)$
with the property that
\[
\{ \tr_{\alpha} a, \tr_{\alpha} b \} = \tr_{\alpha} \lp \overline a,\overline b \rp
\]
for all $a,b\in A$.
\end{thm}

\begin{proof}
Uniqueness is clear since the elements $\tr_{\alpha} a$ generate $\T(A,\alpha)$.

Given $a\in A$, choose a derivation $d_a:A\to A$ inducing $\lp \overline a,-\rp$.
We may assume that $d_a(e_v)=0$ for all $v$. Namely, if $R$ denotes the product
of $m$ copies of $K$ and $\phi:R\to A$ is the $K$-algebra homomorphism sending
$(\lambda_1,\dots,\lambda_m)$ to $\lambda_1 e_1+\dots+\lambda_m e_m$, then $d_a\phi(-)$
is a $K$-derivation $R\to A$, and hence inner since $R$ is a separable $K$-algebra.
Thus $d_a\phi(-) = [a',\phi(-)]$ for some $a'\in A$, and we can replace
the representative $d_a$ by the derivation $x\mapsto d_a(x) - [a',x]$.

By the previous lemma there is a derivation $\psi_a$ of $K[\Rep(A,\alpha)]$
defined by $\psi_a = (d_a)_\alpha$.
Now if $b\in A$ then $\psi_a (b_{ij}) = d_a(b)_{ij}$, so
\[
\psi_a(\tr_{\alpha} b) = \tr_{\alpha} (d_a(b))
= \tr_{\alpha} \lp \overline a,\overline b \rp \in \T(A,\alpha).
\]
For each $f\in \T(A,\alpha)$, choose an expression
\[
f = \sum \lambda^f_{a_1,\dots,a_k} \tr_{\alpha} a_1 \dots \tr_{\alpha} a_k
\]
where the sum is over various collections of elements $a_1,\dots,a_k\in A$,
and the coefficients $\lambda^f_{a_1,\dots,a_k}$ are in $K$.
Since $\psi_a$ is a derivation,
\[
\begin{split}
\psi_a(f)
&= \sum \lambda^f_{a_1,\dots,a_k} \sum_{i=1}^k
\biggl( \prod_{j\neq i} \tr_{\alpha} a_j \biggr) \psi_a(\tr_{\alpha} a_i)
\\
&= \sum \lambda^f_{a_1,\dots,a_k} \sum_{i=1}^k
\biggl( \prod_{j\neq i} \tr_{\alpha} a_j \biggr) \tr_{\alpha} \lp \overline a,\overline a_i \rp,
\end{split}
\]
which shows that $\psi_a$ restricts to a derivation of $\T(A,\alpha)$,
and that this restriction does not depend on the choice of $d_a$.

We define a bracket $\{-,-\}$ on $\T(A,\alpha)$ by
\[
\{ f, g \} = \sum \lambda^f_{a_1,\dots,a_k} \sum_{i=1}^k
\biggl( \prod_{j\neq i} \tr_{\alpha} a_j \biggr) \psi_{a_i} (g).
\]
Clearly $\{f,g\}$ is a derivation in $g$. Moreover,
using that $\lp-,-\rp$ is skew symmetric,
\[
\begin{split}
\{ f, \tr_{\alpha} b \}
&= \sum \lambda^f_{a_1,\dots,a_k} \sum_{i=1}^k
\biggl( \prod_{j\neq i} \tr_{\alpha} a_j \biggr)
\tr_{\alpha} \lp \overline a_i,\overline b\rp
\\
&= - \sum \lambda^f_{a_1,\dots,a_k} \sum_{i=1}^k
\biggl( \prod_{j\neq i} \tr_{\alpha} a_j \biggr)
\tr_{\alpha} \lp\overline b,\overline a_i\rp
\\
&= -\psi_b(f).
\end{split}
\]
Writing
\[
g = \sum \lambda^g_{b_1,\dots,b_\ell} \tr_{\alpha} b_1 \dots \tr_{\alpha} b_\ell,
\]
the fact that $\{f,g\}$ is a derivation in $g$ implies that
\[
\begin{split}
\{f,g\}
&= \sum \lambda^g_{b_1,\dots,b_\ell} \sum_{i=1}^\ell
\biggl( \prod_{j\neq i} \tr_{\alpha} b_j \biggr) \{ f,\tr_{\alpha} b_i\}
\\
&= - \sum \lambda^g_{b_1,\dots,b_\ell}  \sum_{i=1}^\ell
\biggl( \prod_{j\neq i} \tr_{\alpha} b_j \biggr) \psi_{b_i}(f)
\\
&= - \{g,f\}.
\end{split}
\]
As well as showing skew symmetry, this shows that $\{f,g\}$ is a
derivation in $f$, and that it does not depend on the expression for $f$.

Clearly we have
$\{ \tr_{\alpha} a, \tr_{\alpha} b \} = \tr_{\alpha} \lp\overline a,\overline b\rp$.
Since $\pi$ is a noncommutative Poisson structure, we have
\[
\lp \overline a,\lp \overline b,\overline c\rp\rp
+\lp \overline b,\lp \overline c,\overline a\rp\rp
+\lp \overline c,\lp \overline a,\overline b\rp\rp
=0
\]
for all $a,b,c\in A$. This implies that the Jacobi identity
\[
\{ f,\{g,h\}\} + \{ g,\{h,f\}\}+\{ h,\{f,g\}\} = 0
\]
holds for $f = \tr_{\alpha} a$, $g=\tr_{\alpha} b$ and $h=\tr_{\alpha} c$.
Now an induction shows that it holds when $f,g,h$ are products of
elements of the form $\tr_{\alpha} x$, and hence for all $f,g,h\in \T(A,\alpha)$.
Namely, if the Jacobi identity holds for $f_1,g,h$, and for $f_2,g,h$, then
\[
\begin{split}
&\{ f_1 f_2 ,\{g,h\}\} + \{ g,\{h,f_1 f_2\}\}+\{ h,\{f_1 f_2,g\}\}
= f_1 \{ f_2,\{g,h\}\} + \{f_1,\{g,h\}\}f_2
\\
&\quad+ f_1\{g, \{h,f_2\}\} + \{g,f_1\}\{h,f_2\}+\{h,f_1\}\{g,f_2\}+\{g,\{h,f_1\}\}f_2
+ f_1\{h, \{f_2,g\}\}
\\
&\quad+ \{h,f_1\}\{f_2,g\}+ \{f_1,g\}\{h,f_2\}+\{h,\{f_1,g\}\}f_2 =0.
\end{split}
\]
In case $\frac12\notin K$, a similar induction shows that $\{f,f\}=0$
for all $f\in \T(A,\alpha)$.
\end{proof}


\begin{thebibliography}{99}
\bibitem{BG}
J. Block and E. Getzler,
Quantization of foliations, in:
Proceedings of the XXth International Conference on Differential
Geometric Methods in Theoretical Physics, Vol. 1, 2 (New York, 1991), 471--487,
World Sci. Publishing, River Edge, NJ, 1992.

\bibitem{BL}
R. Bocklandt and L. Le Bruyn,
Necklace Lie algebras and noncommutative symplectic geometry,
Math. Z. 240 (2002), 141--167.

\bibitem{CBEG}
W. Crawley-Boevey, P. Etingof and V. Ginzburg,
Noncommutative Geometry and Quiver algebras,
preprint math.AG/0502301.

\bibitem{CBH}
W. Crawley-Boevey and M. P. Holland,
Noncommutative deformations of Kleinian singularities,
Duke Math. J. 92 (1998), 605--635.

\bibitem{FL}
D. R. Farkas and G. Letzter,
Ring theory from symplectic geometry,
J. Pure Appl. Algebra 125 (1998), 155--190.

\bibitem{G}
V. Ginzburg,
Non-commutative symplectic geometry, quiver varieties, and operads,
Math. Res. Lett. 8 (2001), 377--400.

\bibitem{LP}
L. Le Bruyn and C. Procesi,
Semisimple representations of quivers,
Trans. Amer. Math. Soc. 317 (1990), 585--598.

\bibitem{V}
M. Van den Bergh,
Double Poisson algebras,
preprint math.QA/0410528.

\bibitem{Xu}
P. Xu,
Noncommutative Poisson algebras,
Amer. J. Math. 116 (1994), 101--125.
\end{thebibliography}
\end{document}